\documentclass[11pt,a4paper]{article}
\usepackage[T1]{fontenc}
\usepackage[utf8]{inputenc}
\usepackage{lmodern}
\usepackage{amsmath,amssymb,amsthm,mathtools}
\usepackage{microtype}
\usepackage{geometry}
\usepackage{enumitem}
\usepackage{xspace}
\usepackage[colorlinks=true,linkcolor=blue!55!black,citecolor=blue!55!black,urlcolor=blue!55!black]{hyperref}
\usepackage{xcolor}

\usepackage{setspace}

\newif\ifreviewlayout

\ifreviewlayout
\else
\fi	

\newtheorem{theorem}{Theorem}[section]
\newtheorem{proposition}[theorem]{Proposition}
\newtheorem{lemma}[theorem]{Lemma}
\newtheorem{definition}[theorem]{Definition}
\newtheorem{corollary}[theorem]{Corollary}
\theoremstyle{remark}
\newtheorem{example}[theorem]{Example}
\theoremstyle{remark}
\newtheorem{remark}[theorem]{Remark}
\newcommand{\BSTuC}{\hbox{\textnormal{\textsf{BST}}\hspace{-.5pt}\raisebox{.9pt}{$\otimes$}}\xspace}
\newcommand{\BSTC}{\hbox{\textnormal{\textsf{BST}}\hspace{-.7pt}\raisebox{.9pt}{$\times$}}\xspace}

\newcommand{\BST}{\hbox{\textnormal{\textsf{BST}}}\xspace}

\newcommand{\NEXPTIME}{\ensuremath{\mathsf{NEXPTIME}}\xspace}
\newcommand{\NP}{\ensuremath{\mathsf{NP}}\xspace}
\newcommand{\HF}{\ensuremath{\mathsf{HF}}\xspace}
\newcommand{\Pow}{\mathcal{P}}
\newcommand{\Tile}{\mathcal{T}}
\newcommand{\e}{\varnothing}

\title{Cartesian Products in Boolean Set Theory: \NEXPTIME-Completeness of Satisfiability}
\author{Domenico Cantone and
Pietro Ursino \\
\emph{Dipartimento di Matematica e Informatica, Universit\`a di Catania}\\
\emph{Viale Andrea Doria 6, I-95125 Catania, Italy.}  \\
\mbox{E-mail:} \texttt{domenico.cantone@unict.it,
pietro.ursino@unict.it}
}
\date{}
\begin{document}
\maketitle
\begin{abstract}
We prove that the ordinary, finite, and hereditarily finite satisfiability problems for each of the extensions of Boolean set theory obtained by adjoining either the unordered Cartesian product $\otimes$ or the Kuratowski ordered Cartesian product $\times$ are \NEXPTIME-complete. The lower bounds are established by polynomial reductions from exponential-square tiling, based on a binary branching construction admitting hereditarily finite realizations. For the unordered product, these bounds match the upper bounds of Cantone and Ursino. For the ordered product, we establish matching upper bounds, thereby settling the previously open decidability problem.

Our decision procedures for the ordered language rely on finite certificates recording membership profiles and associating sets of admissible target profiles with relevant pairs of input profiles. In conjunction with the local Boolean, disequality, and product conditions, accessibility characterizes ordinary satisfiability, whereas an additional weak-order condition characterizes both finite and hereditarily finite satisfiability. We establish sufficiency by constructing models from the certificates. Since the canonical target relation and the existence of a weak order can be determined from the set of profiles, only this set needs to be guessed for normalized conjunctions. We also show that multivalued targets are indispensable in this characterization and have no uniform bound on their multiplicity, and that ordinary satisfiability admits an optimal uniform bound of $\omega$ on the ranks of the variable values.
\end{abstract}

\noindent\textbf{Keywords:} computable set theory; Boolean set theory; ordered and unordered Cartesian products; ordinary, finite and hereditarily finite satisfiability; \NEXPTIME-completeness; exponential tiling.

\section{Introduction}
Boolean set theory (\BST) is the quantifier-free language of sets with the Boolean operations of union, intersection, and difference. We are concerned with the satisfiability problem for two extensions of this language, obtained by adjoining the unordered or the ordered Cartesian product operator. Given sets $A$ and $B$, their \emph{unordered} and \emph{ordered} Cartesian products are, respectively,
\[
 A\otimes B=\bigl\{\{a,b\}:a\in A,\ b\in B\bigr\},
 \qquad
 A\times B=\bigl\{\langle a,b\rangle:a\in A,\ b\in B\bigr\},
\]
where $\langle a,b\rangle=\bigl\{\{a\},\{a,b\}\bigr\}$ is the Kuratowski ordered pair. Notice that $\{a,b\}=\{b,a\}$, whereas $\langle a,b\rangle=\langle b,a\rangle$ need not hold. 
We denote by $\BSTuC$ and $\BSTC$ the extensions of $\BST$ obtained by adjoining $\otimes$ and $\times$, respectively, and study these two languages separately.

In \cite{CantoneUrsino2024}, Cantone and Ursino established nondeterministic exponential-time upper bounds for ordinary, finite, and hereditarily finite satisfiability of \BSTuC and suggested that their decision procedures might be adapted to the standard ordered Cartesian product, whose satisfiability problem remained open. The complexity of some restricted unordered-product fragments was also investigated by Cantone and Maugeri~\cite{CantoneMaugeri2023}.

In the present paper, we pursue this investigation in two directions. First, we establish the optimality of the known upper bounds for \BSTuC by a polynomial reduction from exponential-square tiling and adapt the reduction to \BSTC. Second, we develop a finite-certificate characterization for \BSTC from which nondeterministic exponential-time decision procedures for the three satisfiability problems follow. More precisely, we prove the following theorem.

\begin{theorem}[Main result]\label{thm:main}
For each of the two languages $\BSTuC$ and $\BSTC$, ordinary, finite, and hereditarily finite satisfiability are \NEXPTIME-complete, even for normalized conjunctions of atomic set constraints and disequalities (in the sense of Lemma~\ref{lem:flatten}).
\end{theorem}
 
To establish the lower bounds, we encode the existence of a tiling of a square of side $2^n$ by a formula of size polynomial in that of the tiling instance. Rather than introducing a variable for each of the $2^{2n}$ grid positions, we associate each position with a binary address of length $2n$. By suitably choosing disjoint sets of tags, we extend each address in both possible ways without identifying elements associated with distinct addresses. Thus every address is represented in every model, even an infinite one, and there is a hereditarily finite realization with exactly one representative per address. To express horizontal and vertical adjacency in the unordered language, we exploit the rule for binary increment: the trailing $1$-bits become $0$, the first subsequent $0$-bit becomes $1$, and the remaining bits are preserved. The corresponding Boolean conditions determine the source and destination of each successor pair without enumerating the grid positions. The reduction for the ordered language uses the same address construction and binary-increment conditions, but the orientation of a pair is already determined by ordered pairing.

To obtain the upper bounds for the ordered product, we associate with each element occurring in a variable value its membership pattern, which we call its \emph{profile}. A formula involves only finitely many variables and hence only finitely many profiles, even if its models are infinite. The Boolean operations can be dealt with at the level of profiles. In treating the Cartesian product, however, one must also account for the fact that pairs whose components have the same respective profiles may themselves have different profiles. We therefore associate a nonempty set of admissible target profiles with each relevant pair of input profiles.

Local profile conditions alone do not ensure the existence of a well-founded model. We require in addition an \emph{accessibility condition}, which allows profiles to be reached from those that do not depend on product generation. Together with the local Boolean, disequality, and product conditions, accessibility suffices for ordinary satisfiability. An additional \emph{weak-order condition} allows the realization to be completed in finitely many steps and yields a hereditarily finite model; thus finite and hereditarily finite satisfiability coincide for the ordered language. Although arbitrarily large target sets may be necessary, we show that a canonical target relation can be computed from the set of profiles and that the existence of a weak order can be tested deterministically. Consequently, only the set of profiles needs to be guessed for normalized conjunctions. The ordinary realization also yields an optimal uniform bound of $\omega$ on the ranks of the variable values.

After introducing the languages and the tiling source problem, Section~\ref{sec:branching} develops the binary address construction. Section~\ref{sec:grid} expresses coordinate equality and binary successor, and Section~\ref{sec:tiling} completes the unordered reduction. Section~\ref{sec:ordered} transfers the lower bound to ordered pairing and establishes the corresponding upper bound through a finite-certificate characterization. We conclude by discussing the scope of the results. Appendices~\ref{app:twobit} and~\ref{app:carry} illustrate the address and successor constructions on small grids.

\section{Languages, semantics, and preliminaries}

\subsection{Syntax and intended semantics}

The atomic formulae of $\BSTuC$ are equalities between terms built from \emph{set variables} $x,y,z,\ldots$ by means of the operations $\cup,\cap,\setminus$ and $\otimes$. The ordered variant $\BSTC$ has the same syntax and semantics, with $\otimes$ replaced by the Kuratowski ordered Cartesian product $\times$. The two symbols designate \emph{alternative} extensions of \BST, not one language containing both operators. In schematic input constraints we shall also use inclusions and abbreviations such as $x\cap y=\e$ and $x\ne\e$, all expressible in the language. There is no membership predicate $\in$ in either input language. Formulas are obtained from atomic formulae by means of the usual propositional connectives.

We generally use lower-case letters for set variables and terms, and upper-case letters for actual sets. In the tiling reduction, however, we also use upper-case letters for schematic set variables, suppressing the assignment notation when referring to their interpretations. We work in the universe $\mathcal{V}$ of well-founded sets, organized into the cumulative hierarchy $(\mathcal{V}_\alpha)_{\alpha\in\mathrm{Ord}}$, where $\mathrm{Ord}$ denotes the class of all ordinals. The hierarchy is defined by $\mathcal{V}_0=\varnothing$, $\mathcal{V}_{\alpha+1}=\mathcal{P}(\mathcal{V}_\alpha)$, and $\mathcal{V}_\lambda=\bigcup_{\alpha<\lambda}\mathcal{V}_\alpha$ for every nonzero limit ordinal $\lambda$. Thus $\mathcal{V}=\bigcup_{\alpha\in\mathrm{Ord}}\mathcal{V}_\alpha$. The rank $\operatorname{rk}(a)$ of a set $a$ is the least ordinal $\alpha$ such that $a\in\mathcal{V}_{\alpha+1}$. In particular, $\HF=\mathcal{V}_\omega=\bigcup_{n<\omega}\mathcal{V}_n$ consists precisely of the \emph{hereditarily finite sets}, namely those sets whose transitive closures are finite, where $\omega$ is the first nonzero limit ordinal, that is, the set of all finite ordinals.

An assignment $M$ maps set variables to sets and extends naturally to all terms. We write $Mt$ for the interpretation of a term $t$ under $M$, using parentheses only when the term is explicitly written as a compound expression. Thus, for example, $M(t\cup u)=Mt\cup Mu$. An assignment is a model of a formula if it satisfies the formula under the usual propositional semantics, with equality interpreted as set equality.
A model is \emph{finite} if each variable occurring in the formula is assigned a finite set, and is \emph{hereditarily finite} (HF) if each such variable is assigned an element of $\HF$. Notice that a finite set need not be hereditarily finite, as is the case with $\{\omega\}$, whereas every element of a hereditarily finite set is itself hereditarily finite. An ordinary model is a model with no finiteness restriction on the values of its variables. Ordinary, finite, and hereditarily finite satisfiability refer to the existence of models of the corresponding kind.

In the semantic proofs, the properties expressed by schematic constraints will be written, for instance, as $A\cap B=\e$ and $A\ne\e$. For example, the formula $x=y\otimes z$ is satisfied by an assignment $M$ precisely when the actual sets $X=Mx$, $Y=My$, and $Z=Mz$ satisfy $X=Y\otimes Z$.

\subsection{Normalization of formulae}
The reduction below is easier to describe using inclusions, intersections, and intermediate set expressions. We measure syntactic size by the number of nodes in a syntax-tree representation, with variable identifiers stored in single machine words, and construction time by unit-cost word operations. Renumbering the variables and encoding their identifiers in binary introduces at most a polynomial overhead in the standard bit model. The following lemma shows that these abbreviations can be eliminated while preserving all three notions of satisfiability.

\begin{lemma}[Flattening]\label{lem:flatten}
Every conjunction of equations, inclusions, and disequalities between set terms constructed from $\cup,\cap,\setminus,\otimes$, with $\e$ allowed as an abbreviation, can be translated into an equisatisfiable conjunction whose conjuncts have the normalized forms
\[
 x=y\cup z,\qquad x=y\setminus z,\qquad x=y\otimes z,
 \qquad x\ne y.
\]
The translation introduces only linearly many conjuncts and fresh variables. No new primitive operation is required. The same statement holds when $\otimes$ is replaced throughout by the ordered Cartesian product $\times$.
\end{lemma}
\begin{proof}
Associate a fresh variable with each nonvariable subterm and add its corresponding defining equation. Eliminate intersections using the set-term identity $x\cap y=x\setminus(x\setminus y)$. Express $x=y$ by $x=y\cup y$, and $x\subseteq y$ by $y=x\cup y$. If the empty set occurs, replace every occurrence of $\e$ by a single fresh variable $z_\e$ constrained by $z_\e=z_\e\setminus z_\e$. Each subterm occurrence requires only a constant number of additional variables and constraints. Hence the resulting conjunction contains only linearly many conjuncts and fresh variables. The equations form a definitional extension. The same translation applies to $\times$. 
All the operations under consideration preserve finiteness and hereditary finiteness. Thus every model of the original conjunction extends to a model of the flattened conjunction in the same semantic class by interpreting each fresh variable as the subterm it represents. Conversely, restricting a model of the flattened conjunction to the original variables gives a model of the original conjunction in the same semantic class.
\end{proof}

We call a conjunction \emph{normalized} if each of its conjuncts has one of the forms listed in Lemma~\ref{lem:flatten}, with the Cartesian product appropriate to the language.

\subsection{The source tiling problem}

We reduce from the origin-constrained exponential-square tiling problem, known to be \NEXPTIME-complete in the precise form used below; see~\cite[Section~4.2, proof of Theorem~17]{BaaderEtAl2016} (and also~\cite{Furer1984,CeylanEtAl2025}). An instance consists of a positive integer $n$ given in unary, a finite nonempty set $D$ of tile types, adjacency relations $H,V\subseteq D\times D$, and a distinguished tile type $d_\star\in D$. Write $\Tile=(D,H,V,d_\star)$ for the tile system and let $L=2^n$.  The question is whether there is a map
\[
 \tau \colon \{0,\ldots,L-1\}^2\longrightarrow D
\]
with $\tau(0,0)=d_\star$ and, for all $x,y\in\{0,\ldots,L-1\}$,
\begin{align}
 x<L-1&\;\Longrightarrow\;(\tau(x,y),\tau(x+1,y))\in H,\label{eq:tileH}\\
 y<L-1&\;\Longrightarrow\;(\tau(x,y),\tau(x,y+1))\in V.\label{eq:tileV}
\end{align}
No wraparound is imposed. The unary encoding of $n$ ensures that the grid has at most exponential size in the length of the input, while allowing our reduction, which uses $2n$ branching stages, to produce a formula of polynomial size.

\section{Finite binary branching and exponential address spaces}
\label{sec:branching}

We begin by showing how to force exponentially many grid representatives while retaining the possibility of a hereditarily finite model. The following lemma relies on disjointness of the parent set and the two tag sets, with no restriction on their cardinalities.

\begin{lemma}[Disjoint-tag branching]\label{lem:branching}
Let $L,T^0,T^1$ be pairwise disjoint sets with $T^0,T^1\ne\e$.
For $A\subseteq L$ and $b\in\{0,1\}$ define
$S_b(A)=A\otimes T^b$. Then, for all $A,B\subseteq L$,
\begin{align}
 S_b(A\cup B)&=S_b(A)\cup S_b(B),\label{eq:branchunion}\\
 S_b(A\cap B)&=S_b(A)\cap S_b(B),\label{eq:branchinter}\\
 S_b(A\setminus B)&=S_b(A)\setminus S_b(B),\label{eq:branchdiff}\\
 S_0(A)\cap S_1(B)&=\e.\label{eq:branchsep}
\end{align}
Moreover, $A\ne\e$ implies $S_b(A)\ne\e$, and $S_b$ is injective.
\end{lemma}
\begin{proof}
For any $a\in L$ and $t\in T^b$, disjointness implies $a\ne t$.
The unordered pair $\{a,t\}$ has exactly one member in $L$, namely
$a$, and exactly one member in $T^b$, namely $t$. Consequently
\[
 \{a,t\}=\{a',t'\},\quad
 a,a'\in L,\quad t,t'\in T^b
 \quad\Longrightarrow\quad a=a',\ t=t'.
\]
The same unordered pair cannot arise from two distinct branches,
because the unique member outside $L$ would then have to belong to
$T^0\cap T^1$. Thus the families of outputs indexed by different
parents or different branch labels are pairwise disjoint, although
a parent may have several children if its tag set has several elements.
Their union descriptions give \eqref{eq:branchunion}--\eqref{eq:branchdiff},
and branch disjointness gives \eqref{eq:branchsep}. Choosing a tag in
$T^b$ proves nonemptiness when $A\ne\e$. Finally, if $A\ne B$, a
parent in their symmetric difference has an output belonging to exactly
one of $S_b(A),S_b(B)$, proving injectivity.
\end{proof}

Fix $k\ge1$. Introduce level variables $L_0,\ldots,L_k$.
For $0\le i<k$ introduce tag variables $T_i^0,T_i^1$, branch
variables $Z_i^0,Z_i^1$, and, for $1\le h\le j\le k$,
bit variables $B_h^j$. Let $\Gamma_k$ contain the constraints
\begin{align}
 L_0&\ne\e,\qquad T_i^0\ne\e,\quad T_i^1\ne\e,
   \label{eq:branchinit}\\
 L_i\cap T_i^b&=\e\quad(b\in\{0,1\}),\qquad
 T_i^0\cap T_i^1=\e,\label{eq:branchdisjoint}\\
 Z_i^b&=L_i\otimes T_i^b\quad(b\in\{0,1\}),\qquad
 L_{i+1}=Z_i^0\cup Z_i^1,\label{eq:branchstage}\\
 B_{i+1}^{i+1}&=Z_i^1,\label{eq:branchnewbit}\\
 B_h^{i+1}&=(B_h^i\otimes T_i^0)\cup
               (B_h^i\otimes T_i^1)
             \quad(1\le h\le i).\label{eq:branchprop}
\end{align}
All constraints indexed by $i$ range over $0\le i<k$ wherever
meaningful; the propagation constraints have $1\le h\le i$.
The displayed intersections, inclusions, empty-set comparisons and
non-flat terms are merely abbreviations, to be eliminated by
Lemma~\ref{lem:flatten}. No cell indexed by a word of length $k$ is
an input variable: the word-indexed cells below are used only in
proofs.

\begin{lemma}[Realization of all addresses]\label{lem:addresses}
In every ordinary model of $\Gamma_k$, for each $0\le i\le k$
there are nonempty, pairwise disjoint cells $\mathcal{L}_\alpha\subseteq L_i$
indexed by $\alpha\in\{0,1\}^i$ such that
\[
 L_i=\mathop{\dot\bigcup}_{\alpha\in\{0,1\}^i}\mathcal{L}_\alpha,
 \qquad
 B_h^i=\bigcup_{\substack{\alpha\in\{0,1\}^i\\\alpha_h=1}}
 \mathcal{L}_\alpha\quad(1\le h\le i).
\]
Here $\epsilon$ denotes the empty word. Thus, writing $G=L_k$ and $B_h=B_h^k$, for every
$\alpha\in\{0,1\}^k$ the region
\[
 \mathcal{R}_\alpha=G\cap
 \bigcap_{\alpha_h=1}B_h\cap
 \bigcap_{\alpha_h=0}(G\setminus B_h)
\]
where an empty intersection is understood relative to $G$, is exactly $\mathcal{L}_\alpha$ and is nonempty. In particular,
$|G|\ge 2^k$ in every model.
\end{lemma}
\begin{proof}
Suppress the assignment notation. Put $\mathcal{L}_\epsilon=L_0$ and,
recursively, $\mathcal{L}_{\alpha b}=\mathcal{L}_\alpha\otimes T_i^b$ whenever
$|\alpha|=i<k$ and $b\in\{0,1\}$. At level $i=0$, the
only cell is nonempty. Suppose the cells at level $i$ are nonempty,
pairwise disjoint and have union $L_i$. By
\eqref{eq:branchinit} and \eqref{eq:branchdisjoint}, Lemma~\ref{lem:branching} applies to
$L_i,T_i^0,T_i^1$. It shows that all the $\mathcal{L}_{\alpha b}$
are nonempty and pairwise disjoint, while
\[
 \bigcup_{\alpha,b}\mathcal{L}_{\alpha b}
 =(L_i\otimes T_i^0)\cup(L_i\otimes T_i^1)
 =L_{i+1}
\]
by \eqref{eq:branchstage} and union preservation. The partition
statement follows by induction, even when the sets are infinite.

At the level at which bit $h$ is created,
\eqref{eq:branchnewbit} gives
$B_h^h=\bigcup_{\alpha\in\{0,1\}^h,\ \alpha_h=1}\mathcal{L}_\alpha$.
If the same identity holds at level $i$, then
\eqref{eq:branchprop} and union preservation give
\[
 B_h^{i+1}=\bigcup_{\substack{\alpha\in\{0,1\}^i\\\alpha_h=1}}
 (\mathcal{L}_{\alpha0}\cup \mathcal{L}_{\alpha1})
 =\bigcup_{\substack{\beta\in\{0,1\}^{i+1}\\\beta_h=1}}
 \mathcal{L}_\beta.
\]
Induction proves the bit identity. The cells form a partition
and their memberships in the final $B_h$ are precisely their
labels, proving the asserted identity for $\mathcal{R}_\alpha$.
\end{proof}

\begin{lemma}[HF realization and size]\label{lem:hfaddresses}
For every $k\ge1$, the conjunction $\Gamma_k$ has a model
whose variable values are hereditarily finite and whose final
cells $\mathcal{L}_\alpha$ are singletons. In this model $|L_k|=2^k$.
Moreover, the conjunction has a polynomial-time normalized
translation of size $O(k^2)$.
\end{lemma}
\begin{proof}
Set $L_0=\{\e\}$. Suppose that $L_i$ and all values
through level $i$ have been constructed as finite subsets of
$\HF$. Choose two distinct elements
$t_i^0,t_i^1\in\HF\setminus L_i$ and set
$T_i^b=\{t_i^b\}$. These tags are nonempty and pairwise
disjoint, and are disjoint from $L_i$. Such choices are always
possible because $L_i$ is finite and $\HF$ is infinite.
Define the branch values, the next level and all propagated bits
by \eqref{eq:branchstage}--\eqref{eq:branchprop}. Finite unions,
differences and products of finite hereditarily finite sets remain
hereditarily finite, so every value lies in $\HF$.
By Lemma~\ref{lem:branching}, singleton tags give exactly one child per parent in each branch, and distinct parent--branch choices give distinct children. Thus the cells remain singletons and after $k$ steps their number is $2^k$.

There are $O(k)$ level, tag and branch variables and
$\sum_{i=1}^k i=O(k^2)$ bit variables. There are $O(k^2)$
propagation operations; the remaining constraints have linear
total size. Lemma~\ref{lem:flatten} yields the normalized form
without asymptotic increase, and its defining extensions
preserve hereditary finiteness (as well as finite and ordinary
satisfiability).
\end{proof}

\section{Encoding a grid by Boolean profiles}
\label{sec:grid}
Fix a tiling instance $(n,D,H,V,d_\star)$ and put $k=2n$.
Adjoin $\Gamma_k$ to the formula under construction. The final level
and bit variables are
\[
 G=L_k,\qquad B_j=B_j^k,\qquad
 \bar B_j=G\setminus B_j\quad(1\le j\le 2n).
\]
Every $p\in G$ has an address profile
\[
 \beta(p)=(\beta_1(p),\ldots,\beta_{2n}(p))\in\{0,1\}^{2n},
 \quad\beta_j(p)=1\ \Longleftrightarrow\ p\in B_j.
\]
Write $\beta(p)=(\beta^h(p),\beta^v(p))$, where the first and last
$n$ bits encode respectively horizontal and vertical coordinates,
least significant bit first. Thus
\[
 x(p)=\sum_{j=1}^{n}2^{j-1}\beta_j(p),\qquad
 y(p)=\sum_{j=1}^{n}2^{j-1}\beta_{n+j}(p).
\]
By Lemma~\ref{lem:addresses}, every address profile is realized in
$G$ (regardless of whether $G$ is finite). Consequently
\begin{equation}\label{eq:surject}
 \{(x(p),y(p)):p\in G\}=\{0,\ldots,2^n-1\}^2.
\end{equation}
Although an address may have several representatives, the following sets allow us to compare their bits independently of the choice of representatives.

For $1\le j\le 2n$ and $J\subseteq\{1,\ldots,2n\}$, define
\begin{equation}\label{eq:Qdef}
 Q_j=(B_j\otimes B_j)\cup(\bar B_j\otimes\bar B_j),
 \qquad Q_J=\bigcap_{j\in J}Q_j.
\end{equation}
For $J=\e$, put $Q_J=G\otimes G$; this convention keeps all
definitions within the permitted set operations.

\begin{lemma}[Equality of bits]\label{lem:eqbits}
For any $p,q\in G$ and any $J\subseteq\{1,\ldots,2n\}$,
\begin{equation}\label{eq:iffbits}
 \{p,q\}\in Q_J
 \quad\Longleftrightarrow\quad
 \beta_j(p)=\beta_j(q)\quad\text{for every }j\in J.
\end{equation}
\end{lemma}
\begin{proof}
For a given $j$, the pair $\{p,q\}$ belongs to
$B_j\otimes B_j$ exactly when both entries are in $B_j$, and to
$\bar B_j\otimes\bar B_j$ exactly when both entries are outside $B_j$
but inside $G$.  This is equivalent to equality of the $j$th bits. For $p=q$ the
unordered pair becomes the singleton $\{p\}$, which is in exactly
one of these two products, according to the bit value. Intersecting
the conditions proves the assertion for every $J$, including
$J=\varnothing$ under the stated convention.
\end{proof}

To express adjacency, consider a possible \emph{carry pivot} $i$ ($1\le i\le2n$), namely the position of the first zero bit in the coordinate to be incremented, read from the least significant bit. Let
\[
 I_i=\begin{cases}
 \{1,\ldots,i-1\},&1\le i\le n,\\
 \{n+1,\ldots,i-1\},&n<i\le 2n,
 \end{cases}
\]
and define source and destination sets
\begin{align}
 S_i&=\bar B_i\cap\bigcap_{j\in I_i}B_j,\label{eq:Sdef}\\
 D_i&=B_i\cap\bigcap_{j\in I_i}\bar B_j.\label{eq:Ddef}
\end{align}
An empty intersection in \eqref{eq:Sdef}--\eqref{eq:Ddef} is
interpreted as $G$.  Define the indices of the preserved bits by
\[
 J_i=\begin{cases}
 \{i+1,\ldots,n\}\cup\{n+1,\ldots,2n\},&1\le i\le n,\\
 \{1,\ldots,n\}\cup\{i+1,\ldots,2n\},&n<i\le 2n.
 \end{cases}
\]
Finally put
\begin{equation}\label{eq:Rdef}
 R_i=(S_i\otimes D_i)\cap Q_{J_i}.
\end{equation}
Note that $S_i\cap D_i=\e$.

\begin{lemma}[Directed successor from unordered pairs]\label{lem:successor}
Let $p,q\in G$ and $1\le i\le 2n$. Then $\{p,q\}\in R_i$
if and only if there is an ordering $(s,t)$ of $p,q$ such that
\[
 s\in S_i,\qquad t\in D_i,\qquad
 \begin{cases}
 x(t)=x(s)+1,\quad y(t)=y(s),&1\le i\le n,\\
 x(t)=x(s),\quad y(t)=y(s)+1,&n<i\le 2n.
 \end{cases}
\]
Whenever such an ordering exists, it is unique. Moreover,
$\beta_i(s)=0$ and $\beta_i(t)=1$, so $i$ is the carry
position of the corresponding binary increment. Thus the
source and destination are determined unambiguously even
though $\otimes$ forms unordered pairs.
\end{lemma}
\begin{proof}
First suppose $1\le i\le n$. Since $S_i\cap D_i=\e$,
membership in $S_i\otimes D_i$ is equivalent to the existence
of a unique orientation $(s,t)$ with $s\in S_i$ and $t\in D_i$.
In this orientation, the horizontal bits below $i$ are $1$ at
$s$ and $0$ at $t$, while the pivot bit is $0$ at $s$ and $1$
at $t$. By Lemma~\ref{lem:eqbits}, membership in $Q_{J_i}$
means that all higher horizontal bits and all vertical bits
are equal at $s$ and $t$. Consequently,
\[
 x(t)-x(s)=2^{i-1}-\sum_{j=1}^{i-1}2^{j-1}=1,
 \qquad y(t)=y(s).
\]
Conversely, suppose $s\in S_i$, $t\in D_i$,
$x(t)=x(s)+1$, and $y(t)=y(s)$. The bits below and at $i$
have the prescribed values. Since the total difference of
the horizontal coordinates is $1$, their higher horizontal
bits must agree by uniqueness of binary expansion; equality
of vertical coordinates likewise forces all vertical bits to
agree. Thus $\{s,t\}\in Q_{J_i}$, and therefore
$\{p,q\}=\{s,t\}\in R_i$.
The proof for $n<i\le 2n$ is identical, exchanging the
horizontal and vertical coordinates and using the lower
vertical bits $n+1,\ldots,i-1$.
\end{proof}

\section{Tiling constraints and proof of \NEXPTIME-hardness}
\label{sec:tiling}
For every $d\in D$ introduce a variable $C_d$ for the set of grid
representatives assigned tile $d$.  Require
\begin{equation}\label{eq:partition}
 G=\bigcup_{d\in D}C_d,\qquad
 C_d\cap C_e=\e\quad(d\ne e).
\end{equation}
These conditions in particular imply $C_d\subseteq G$ for every $d$.
Put $O=G\setminus\bigcup_{j=1}^{2n}B_j$. By \eqref{eq:surject}, $O$ is nonempty, and it consists precisely of
representatives of the grid origin.  Add the origin constraint
\begin{equation}\label{eq:origin}
 O\subseteq C_{d_\star}.
\end{equation}

To rule out forbidden adjacencies, we restrict the successor-pair sets $R_i$ to pairs whose source carries tile $d$ and whose destination carries tile $e$. These restrictions can be expressed without introducing separate relation variables.
For each $1\le i\le n$ and each forbidden horizontal adjacency
$(d,e)\in D^2\setminus H$ impose
\begin{equation}\label{eq:noH}
 \bigl((S_i\cap C_d)\otimes(D_i\cap C_e)\bigr)
 \cap Q_{J_i}=\e.
\end{equation}
For each $n<i\le 2n$ and $(d,e)\in D^2\setminus V$ impose
\begin{equation}\label{eq:noV}
 \bigl((S_i\cap C_d)\otimes(D_i\cap C_e)\bigr)
 \cap Q_{J_i}=\e.
\end{equation}
Let $\Psi_{n,\Tile}$ be the conjunction of $\Gamma_{2n}$,
\eqref{eq:partition}, \eqref{eq:origin}, and all instances of
\eqref{eq:noH}--\eqref{eq:noV}, together with the definitional equations
for the auxiliary set terms.

\begin{lemma}[Soundness of the reduction]\label{lem:sound}
If $\Psi_{n,\Tile}$ is satisfiable, the given exponential-square tiling
instance has a solution.
\end{lemma}
\begin{proof}
Take any model of $\Psi_{n,\Tile}$.  By \eqref{eq:surject}, for
every $(x,y)\in\{0,\ldots,2^n-1\}^2$ we may choose a representative
$p_{x,y}\in G$ of these coordinates.  By \eqref{eq:partition}, exactly
one $d\in D$ has $p_{x,y}\in C_d$.  Define $\tau(x,y)=d$.
All representatives of $(0,0)$ belong to $O$, hence
$\tau(0,0)=d_\star$ by \eqref{eq:origin}.

Suppose $x<2^n-1$.  Let $p=p_{x,y}$ and $q=p_{x+1,y}$.
There is a unique pivot $i\in\{1,\ldots,n\}$ for binary increment
from $x$ to $x+1$. By the definitions of $S_i$ and $D_i$,
$p\in S_i$ and $q\in D_i$. Moreover, the bits at all positions
in $J_i$ agree, so Lemma~\ref{lem:eqbits} gives $\{p,q\}\in Q_{J_i}$.
If $(\tau(x,y),\tau(x+1,y))=(d,e)\notin H$, this unordered pair
belongs to the left-hand side of \eqref{eq:noH}, a contradiction.
Thus the horizontal constraint \eqref{eq:tileH} is satisfied.
Proceeding analogously with a vertical pivot, we obtain \eqref{eq:tileV}. Thus $\tau$ is a valid tiling.
\end{proof}

\begin{lemma}[HF completeness of the reduction]\label{lem:complete}
If the exponential-square tiling instance has a solution, then
$\Psi_{n,\Tile}$ has a hereditarily finite model and therefore
also a finite and an ordinary model.
\end{lemma}
\begin{proof}
Let $\tau$ be a valid tiling and take the hereditarily finite
model of $\Gamma_{2n}$ supplied by Lemma~\ref{lem:hfaddresses}.
In this model every final address cell has exactly one element.
Interpret $G,B_j,\bar B_j,Q_j,S_i,D_i$ and $R_i$ by their
set-theoretic definitions. Every $p\in G$ has well-defined
coordinates $x(p),y(p)$. For every $d\in D$ put
\[
 C_d=\{p\in G:\tau(x(p),y(p))=d\}.
\]
These are pairwise disjoint hereditarily finite sets whose union is
$G$, so \eqref{eq:partition} holds; the origin condition
\eqref{eq:origin} also holds because $\tau(0,0)=d_\star$.
Suppose some forbidden horizontal intersection
\eqref{eq:noH} were nonempty. It would contain an unordered pair
$\{p,q\}$ with
$p\in S_i\cap C_d$, $q\in D_i\cap C_e$, and
$\{p,q\}\in Q_{J_i}$. Since $S_i\cap D_i=\e$, the
orientation is unique. Lemma~\ref{lem:successor} implies
\[
 (x(q),y(q))=(x(p)+1,y(p)),
\]
so the definition of the colors and the horizontal tiling rule
force $(d,e)\in H$, contradicting the choice of the forbidden
pair. The vertical clauses follow identically.
All auxiliary variables can be interpreted as the values of
their defining set terms. Since finite Boolean operations and
unordered products preserve $\HF$, all such values are
hereditarily finite. Thus $\Psi_{n,\Tile}$ has an HF model.
\end{proof}

\begin{lemma}[Size of the reduction]\label{lem:size}
There is a deterministic polynomial-time algorithm which, given
$(n,D,H,V,d_\star)$ with $n$ in unary, constructs a normalized
$\BSTuC$ conjunction equivalent to $\Psi_{n,\Tile}$ with respect
to each of the three notions of satisfiability.
\end{lemma}
\begin{proof}
The finite binary-address part $\Gamma_{2n}$ consists of
$O(n^2)$ constraints.  There are $2n$ bit sets, $2n$ equal-bit sets $Q_j$,
and $2n$ source and destination pairs.  Each of these auxiliary terms
requires $O(n)$ operations or fewer.  The disjoint-cover constraints for $G$ require
$O(|D|^2)$ disjointness conditions.  There are at most
$2n|D|^2$ forbidden adjacency constraints; each is an expression
of size $O(n)$, or shorter if previously introduced terms are shared.
The level and bit variables are linked by defining equations: their
recursively described values are not expanded into a single term.
Summing the above contributions, even without sharing repeated subexpressions, we obtain a total length of $O(n^2|D|^2+n^2+|D|^2)$, and the formula can be constructed
in polynomial time.  Lemma~\ref{lem:flatten} puts it into the
specified flat syntax without changing satisfiability in any of
the three semantics.
\end{proof}

\begin{proof}[Unordered-product case of Theorem~\ref{thm:main}]
By Lemmas~\ref{lem:sound}--\ref{lem:size}, there is a polynomial-time
many-one reduction from origin-constrained exponential-square tiling
into normalized $\BSTuC$ conjunctions. The soundness lemma applies to
every ordinary model, and the completeness lemma constructs an HF
model from every positive tiling instance. Therefore, for each of the
three semantics separately, the normalized satisfiability problem is
\NEXPTIME-hard. Cantone and Ursino~\cite{CantoneUrsino2024}
proved nondeterministic exponential-time upper bounds for ordinary
and (hereditarily) finite $\BSTuC$ satisfiability. These bounds also
apply to normalized conjunctions, which form a sublanguage. Hence all
three normalized satisfiability problems are \NEXPTIME-complete. The
same lower bounds apply to the full quantifier-free language.
The ordinary upper bound for arbitrary formulae is established in the
cited work; the finite and HF upper bounds extend from conjunctions
by the usual nondeterministic reduction to conjunctions of literals.
\end{proof}

\section{The ordered Cartesian product: lower and upper bounds}
\label{sec:ordered}

We turn now to the language $\BSTC$. Two properties of Kuratowski pairing will be used in what follows: its injectivity and the fact that the rank of a pair strictly exceeds the ranks of both components. More precisely,
\begin{equation}\label{eq:kuratowski}
 \langle a,b\rangle=\langle c,d\rangle
 \Longleftrightarrow (a=c\ \text{and}\ b=d),\qquad
 \operatorname{rk}\langle a,b\rangle
 =\max\{\operatorname{rk}a,\operatorname{rk}b\}+2.
\end{equation}
Every ordered pair has cardinality at most two. The rank formula
follows since $\operatorname{rk}\{a\}=\operatorname{rk}(a)+1$
and $\operatorname{rk}\{a,b\}=\max\{\operatorname{rk}(a),
\operatorname{rk}(b)\}+1$. We use these properties in the
construction of finite and infinite models.

\subsection{Transferring the exponential tiling reduction}

\begin{lemma}[Ordered binary branching]\label{lem:orderedbranch}
Let $T^0,T^1$ be disjoint nonempty sets. For $b\in\{0,1\}$, define
$F_b(A)=A\times T^b$. For arbitrary sets $A,B$,
\[
 \begin{aligned}
 F_b(A\cup B)&=F_b(A)\cup F_b(B),&
 F_b(A\cap B)&=F_b(A)\cap F_b(B),\\
 F_b(A\setminus B)&=F_b(A)\setminus F_b(B),&
 F_0(A)\cap F_1(B)&=\varnothing.
 \end{aligned}
\]
Moreover, $F_b$ is injective, and it maps nonempty sets to nonempty sets.
No disjointness between $A,B$ and the tags is necessary.
\end{lemma}
\begin{proof}
Injectivity of ordered pairing implies that the fibers
$\{a\}\times T^b$ for distinct $a$ are disjoint, and that the two
collections of fibers for $b=0,1$ are disjoint.  The stated identities
follow by taking unions of fibers.  Nonemptiness follows from
$T^b\ne\varnothing$, and if $A\ne B$ a fiber indexed by an element
of the symmetric difference distinguishes $F_b(A)$ and $F_b(B)$.
\end{proof}

\begin{proposition}[Ordered-product lower bound]
\label{prop:orderedhard}
Ordinary, finite, and hereditarily finite satisfiability of
$\BSTC$ are \NEXPTIME-hard, even for normalized
conjunctions.
\end{proposition}
\begin{proof}
In view of the injectivity of ordered pairing, there is no need to recover the orientation of a pair from its components. We translate the reduction of Sections~\ref{sec:branching}--\ref{sec:tiling}, in polynomial time, by replacing every occurrence of $\otimes$ in the branching constraints and in the grid encoding by $\times$.
Lemma~\ref{lem:orderedbranch} gives the counterparts of
Lemmas~\ref{lem:addresses} and~\ref{lem:hfaddresses}: for $k=2n$
all $2^{2n}$ Boolean address cells are forced to be nonempty in
\emph{every} model, and the construction has an HF model with
one element per cell.  The tags can be chosen as distinct HF
singletons, as in Lemma~\ref{lem:hfaddresses}.

More explicitly, with the same final bit sets $B_j$ and
$\bar B_j=G\setminus B_j$, for $1\le j\le 2n$ and $J\subseteq\{1,\ldots,2n\}$ set
\[
 Q_j^{\times}=(B_j\times B_j)\cup
                (\bar B_j\times\bar B_j),\qquad
 Q_J^{\times}=\bigcap_{j\in J}Q_j^{\times}.
\]
As before, the empty intersection is interpreted as the ambient product $G\times G$. For $p,q\in G$, $\langle p,q\rangle\in Q_J^{\times}$ exactly when
$p$ and $q$ agree on all bit positions in $J$. Using the source and
destination sets $S_i,D_i$ from \eqref{eq:Sdef}--\eqref{eq:Ddef}, put
\[
 R_i^{\times}=(S_i\times D_i)\cap Q_{J_i}^{\times}.
\]
The binary-increment calculation in Lemma~\ref{lem:successor}
shows that an ordered pair $\langle p,q\rangle$ belongs to
$R_i^{\times}$ precisely when $p$ and $q$ encode an
$i$-pivot horizontal or vertical successor, with $p$ the source.
For every forbidden horizontal pair $(d,e)\in D^2\setminus H$ add
\[
 ((S_i\cap C_d)\times(D_i\cap C_e))\cap Q_{J_i}^{\times}
 =\varnothing \qquad (1\le i\le n),
\]
and impose the same constraints for $(d,e)\in D^2\setminus V$ and $n<i\le 2n$.
The disjoint-cover and origin constraints are unchanged.

If these constraints admit any model, choosing a representative
of each address yields a tiling: every forbidden adjacent pair
would contradict one of the displayed constraints. Conversely,
any legal tiling colors the singleton address cells of the HF
branching model and satisfies every constraint.  The number of
stages, bit variables, and forbidden-pair clauses remains
polynomial in the tiling instance.  Flattening replaces each
nested term by a fresh variable, preserving all three semantics.
Thus we obtain the three claimed polynomial-time reductions.
\end{proof}

\subsection{Finite certificates for ordered-product satisfiability}

Inspired by the profile-graph method of Cantone and Ursino~\cite{CantoneUrsino2024}, we characterize satisfiability for the ordered Cartesian product in terms of finite certificates. We shall first formulate the conditions on profiles and then show how they allow us to construct ordinary or hereditarily finite models. On the grounds of this characterization, we shall derive the decision procedures.

Fix a normalized conjunction $\Phi$ with variable set $W$,
$|W|=m$, whose literals have the forms
\[
 x=y\cup z,\quad x=y\setminus z,\quad x=y\times z,
 \quad x\ne y.
\]
Write $\mathcal M$ for the set of product literals $x=y\times z$.
A \emph{profile} is a nonzero Boolean map $p:W\to\{0,1\}$. Given an assignment $M$, its \emph{profile block} associated with $p$ is
\[
 B_p^M=\{a\in\bigcup_{x\in W}Mx:\ \forall x\in W\ (a\in Mx\iff p(x)=1)\}.
\]
The nonempty profile blocks partition $\bigcup_{x\in W}Mx$. In the abstract realization constructions below, $B_p$ denotes the block being constructed for $p$.
In constructing a model, we distinguish the profiles whose elements are not required to belong to any product-variable value: their blocks can be initialized with elements that are not ordered pairs. Accordingly, for a set $P$ of profiles, define its \emph{source profiles} by
\begin{equation}\label{eq:sourceprofiles}
 S(P)=\{p\in P:p(x)=0\text{ for every }(x=y\times z)\in\mathcal M\}.
\end{equation}
On the other hand, pairs formed from elements in two blocks must be assigned to target blocks whenever a product equation requires them to belong to its left-hand side. For $p,q\in P$ call $(p,q)$ \emph{active} if
$p(y)=q(z)=1$ for at least one product literal $x=y\times z$.
Let $N(P)$ be the set of active ordered pairs of profiles, also called \emph{active nodes}.
An inactive input-profile pair cannot contribute an element to any
product-variable value, and therefore needs no target set.

\begin{definition}[Ordered profile certificate]\label{def:cert}
An \emph{ordered profile certificate} for $\Phi$ consists of
a (possibly empty) set
$P\subseteq\{0,1\}^W\setminus\{0\}$ and, for each
$(p,q)\in N(P)$, a nonempty set $T(p,q)\subseteq P$ of
admissible target profiles, satisfying the following conditions:
\begin{enumerate}[label=(\roman*)]
\item \label{cert:boolean}For every $p\in P$, each Boolean equation
of $\Phi$ holds bitwise:
\[
 p(x)=p(y)\lor p(z)\quad(x=y\cup z),\qquad
 p(x)=p(y)\land\neg p(z)\quad(x=y\setminus z).
\]
For every disequality $x\ne y$ there is a $p\in P$ with
$p(x)\ne p(y)$.
\item \label{cert:products}For each active $(p,q)$ and each
$r\in T(p,q)$, \emph{all} product equations of $\Phi$ obey
\begin{equation}\label{eq:targetprofile}
 r(x)=p(y)q(z)\quad\text{for every }(x=y\times z)\in\mathcal M.
\end{equation}
\item \label{cert:access}A subset $C\subseteq P$ is closed if
\[
 p,q\in C,\quad (p,q)\in N(P)
 \quad\Longrightarrow\quad T(p,q)\subseteq C.
\]
The smallest closed subset of $P$ containing $S(P)$ must be $P$ itself. We call this condition \emph{accessibility}.
\end{enumerate}
The certificate is \emph{weakly ordered} if, additionally, there
is a total order $<$ on $P$ such that
\begin{equation}\label{eq:weakorder}
 \max_{<}\{p,q\}<\max_{<}T(p,q)
 \qquad\text{for every }(p,q)\in N(P).
\end{equation}
Such a total order is called a \emph{weak order} for the certificate.
For $P=\varnothing$, conditions (i)--(iii) hold vacuously
except that the absence of disequalities is necessary. This is the
degenerate certificate of the all-empty assignment: every normalized
Boolean and product equation becomes $\varnothing=\varnothing$.
\end{definition}

\begin{remark}[Admissible target profiles]
The target sets $T(p,q)$ specify admissible destinations, not necessarily all destinations that must be realized. Every pair of elements from an active input node must lie in an admissible target block, and a target block can receive pairs from several nodes. Condition~\ref{cert:products} constrains all product literals simultaneously. Example~\ref{ex:multitarget} illustrates why a single-valued map would not in general record all targets realized by a given model; Proposition~\ref{prop:targetmultiplicity} shows that even existence of a certificate may require multiple targets.

Notice that, by condition~\ref{cert:products}, for every active $(p,q)$, some product literal $x=y\times z$ has
$p(y)q(z)=1$, and hence every $r\in T(p,q)$ satisfies $r(x)=1$.
Thus
\begin{equation}\label{eq:targets-not-sources}
 T(p,q)\cap S(P)=\varnothing\qquad((p,q)\in N(P)).
\end{equation}
In particular, no product-generation step adds an element to a source
profile block.
\end{remark}

\begin{example}[A simple finite certificate]\label{ex:simplecert}
For the conjunction $x=y\times z$, set $a=3$, $b=4$ and
$c=\langle a,b\rangle$, with $3,4$ finite von Neumann ordinals.
Take $My=\{a\}$, $Mz=\{b\}$, and $Mx=\{c\}$.
Writing profiles in the coordinate order $(x,y,z)$, we obtain
\[
 \alpha=(0,1,0),\quad \beta=(0,0,1),\quad
 \gamma=(1,0,0),\qquad
 P=\{\alpha,\beta,\gamma\}.
\]
The sources are $S(P)=\{\alpha,\beta\}$, and the only active
node is $(\alpha,\beta)$, with
$T(\alpha,\beta)=\{\gamma\}$. It is accessible in one step,
and the order $\alpha<\beta<\gamma$ satisfies the weak-order condition.
This elementary example separates the \emph{Boolean data} (the three
profiles), the \emph{product data} (a single active node), and the
global accessibility and order requirements.
\end{example}

\begin{example}[Multiple targets and weak ordering]
\label{ex:multitarget}
Consider the single equation $x=y\times z$. Let $a=3$ and $b=4$
be finite von Neumann ordinals and define
\[
 u_1=\langle a,b\rangle,\qquad
 u_2=\langle u_1,b\rangle,\qquad
 u_3=\langle u_2,b\rangle.
\]
The assignment
\[
 My=\{a,u_1,u_2\},\qquad
 Mz=\{b\},\qquad
 Mx=\{u_1,u_2,u_3\}
\]
is an HF model. Its four nonzero profiles, in coordinate order
$(x,y,z)$, are
\[
\begin{array}{c|ccc|l}
 &x&y&z&\text{elements}\\\hline
 \alpha&0&1&0&a\\
 \beta &0&0&1&b\\
 \gamma&1&1&0&u_1,u_2\\
 \delta&1&0&0&u_3
\end{array}
\]
The sources are $\alpha,\beta$, and the only active nodes are $(\alpha,\beta)$ and $(\gamma,\beta)$, with
\[
 T(\alpha,\beta)=\{\gamma\},\qquad
 T(\gamma,\beta)=\{\gamma,\delta\}.
\]
The second target set really has two elements: pairing $u_1$ with
$b$ produces $u_2$ of profile $\gamma$, whereas pairing $u_2$
with $b$ produces $u_3$ of profile $\delta$.
The order $\alpha<\beta<\gamma<\delta$ satisfies the weak-order
condition: the maximal target of $(\gamma,\beta)$ is $\delta>\gamma$.
But no order can make \emph{every} target of this node larger
than its inputs, because $\gamma$ is itself a target.
\end{example}

\begin{proposition}[Indispensability of multivalued targets]\label{prop:targetmultiplicity}
For every integer $k\ge 2$, there is an HF-satisfiable normalized $\BSTC$ conjunction for which every ordered profile certificate has a node with at least $k$ targets.
\end{proposition}
\begin{proof}
For $k=2$, consider the normalized conjunction
\[
\begin{gathered}
 x=y\times z,\qquad x=u\cup v,\qquad u=u\setminus v,\\
 y=y\setminus x,\qquad z=z\setminus x,\qquad y=y\setminus z,\\
 e=e\setminus e,\qquad u\ne e,\qquad v\ne e.
\end{gathered}
\]
An HF model is given by $My=\{3,4\}$, $Mz=\{5\}$, $Mu=\{\langle 3,5\rangle\}$, $Mv=\{\langle 4,5\rangle\}$, $Mx=Mu\cup Mv$, and $Me=\e$, where the numbers denote finite von Neumann ordinals. The two ordered pairs are distinct by injectivity and neither is one of the three ordinals, since a Kuratowski pair has at most two elements.

The Boolean equations allow exactly the following nonzero profiles:
\[
\begin{array}{c|cccccc}
 &x&y&z&u&v&e\\\hline
 p&0&1&0&0&0&0\\
 q&0&0&1&0&0&0\\
 r&1&0&0&1&0&0\\
 s&1&0&0&0&1&0
\end{array}
\]
The disequalities force $r,s\in P$. Neither is a source, and the only possible active node is $(p,q)$. Accessibility therefore forces $p,q\in P$ and $r,s\in T(p,q)$. Condition~\ref{cert:products} excludes the source profiles from $T(p,q)$, so $T(p,q)=\{r,s\}$.

For arbitrary $k\ge 2$, replace $u,v$ by $u_1,\ldots,u_k$ and impose $u_i\ne e$ for every $i$ and $u_i=u_i\setminus u_j$ for $i<j$. Retain the product equation, the three disjointness equations for $x,y,z$, and $e=e\setminus e$. Express $x=\bigcup_{i=1}^k u_i$ by introducing partial-union variables $a_2,\ldots,a_{k-1}$, with $a_1$ denoting $u_1$ and $a_k$ denoting $x$, and imposing $a_i=a_{i-1}\cup u_i$ for $2\le i\le k$. This is a normalized conjunction. Its HF satisfiability follows by taking $My$ to consist of $k$ distinct finite ordinals of cardinality at least three, $Mz$ to consist of a further such ordinal, and each $Mu_i$ to consist of a distinct pair from $My\times Mz$. Interpret the partial unions accordingly.

The Boolean constraints now allow precisely two input profiles (membership in $y$ alone or in $z$ alone) and $k$ output profiles (membership in $x$ and exactly one $u_i$, with all partial-union bits determined by that choice). Every output profile must occur by the disequalities. None is a source, and the only possible active node has the two input profiles as its inputs. Accessibility forces both inputs and all $k$ output profiles to occur, with all the latter in the target set of that node.
\end{proof}

The proposition concerns certificates based on the membership profiles of the given variables; it does not rule out different representations using more refined combinatorial data. Nor does it conflict with choosing a single designated target per node during completion: the initialization may need several targets.

\begin{example}[Necessity of accessibility]
\label{ex:noaccess}
Consider the normalized conjunction $x=x\times x$, $z=z\setminus z$, and $x\ne z$.
If accessibility were omitted, the certificate with $P=\{p\}$,
$p(x)=1$, $p(z)=0$, and $T(p,p)=\{p\}$ would satisfy
all remaining local profile and target conditions: it would
merely describe a circular dependence of $p$ on itself.
To see that no model exists, suppose by way of contradiction that $M$ is a model. Since $Mx$ is nonempty, choose
an element $a\in Mx$ of minimal rank. Since $Mx=Mx\times Mx$, we
would have $a=\langle b,c\rangle$ for $b,c\in Mx$, both of
strictly smaller rank than $a$, a contradiction. The spurious
certificate is rejected because $S(P)=\varnothing$, so the
accessibility closure cannot reach $p$.
\end{example}

\begin{example}[Ordinary satisfiability need not be finite]
\label{ex:ordinaryonly}
Consider the normalized conjunction
\begin{equation}\label{eq:ordinaryonly}
 y=w\cup x,\qquad x=y\times z,\qquad
 w=w\setminus x,\qquad v=v\setminus v,\qquad w\ne v,\qquad z\ne v.
\end{equation}
Choose $a=3$, $b=4$ and define $u_0=a$ and
$u_{i+1}=\langle u_i,b\rangle$ for $i\ge0$. The assignment
\[
 Mw=\{a\},\qquad Mz=\{b\},\qquad
 Mx=\{u_i:i\ge1\},\qquad
 My=\{u_i:i\ge0\},\qquad Mv=\e
\]
is an ordinary well-founded model. Its nonzero profiles in the
coordinate order $(x,y,z,w,v)$ are
\[
 \alpha=(0,1,0,1,0),\qquad
 \beta=(0,0,1,0,0),\qquad
 \gamma=(1,1,0,0,0).
\]
The source profiles are $\alpha,\beta$, while the two active nodes have target sets
$T(\alpha,\beta)=\{\gamma\}$ and
$T(\gamma,\beta)=\{\gamma\}$. Thus the certificate is
accessible, but it has no weak order: the second node would require
$\max_<\{\gamma,\beta\}<\gamma$, which is impossible in a
total order.

Indeed, no finite model exists. In any finite model $M$, the sets $My$ and $Mz$
are nonempty. Choose $c\in My$ of maximal rank and any $d\in Mz$.
Then $\langle c,d\rangle\in My\times Mz=Mx\subseteq My$, while
$\operatorname{rk}\langle c,d\rangle>\operatorname{rk}(c)$,
contradicting maximality.
\end{example}

We next establish the characterization suggested by the preceding examples, starting with the conditions induced on profiles by a model.

\begin{lemma}[Necessity of the certificate]\label{lem:certnecessary}
Every ordinary model of $\Phi$ induces an ordered profile
certificate.  If the model is finite, the certificate can be chosen
weakly ordered.
\end{lemma}
\begin{proof}
Let $M$ be a model.  For each element $a$ of
$\bigcup_{x\in W}Mx$ let
$p_a(x)=1$ exactly when $a\in Mx$, and let
$P=\{p_a:a\in\bigcup_{x\in W}Mx\}$.  This is a set of nonzero profiles
of cardinality at most $2^m-1$.  Boolean equations and disequalities
give condition~\ref{cert:boolean}.
For each active $(p,q)$ define $T(p,q)$ to be the set of profiles
of $\langle a,b\rangle$ as $a,b$ range over elements of respective
profiles $p,q$.  Activity guarantees that every such pair belongs
to at least one product-variable value, so its profile is in $P$;
therefore $T(p,q)$ is nonempty.  Equality
$Mx=My\times Mz$, together with injectivity of ordered
pairing, yields \eqref{eq:targetprofile} for every product literal,
including those for which $p(y)q(z)=0$.

We verify accessibility by well-founded induction on the rank of
an element $a\in\bigcup_{x\in W}Mx$.  If $p_a\in S(P)$, its profile is
present from the start.  Otherwise $a\in Mx$ for at least one
product literal $x=y\times z$, and hence
$a=\langle b,c\rangle$ with $b\in My$ and $c\in Mz$.
By \eqref{eq:kuratowski}, both $b,c$ have rank strictly below the
rank of $a$.  Their profiles are accessible by induction, the
ordered pair of these profiles is active, and the definition of
$T$ puts $p_a$ among its targets.  This proves~\ref{cert:access}. In particular, if $P$ is nonempty,
then $S(P)$ is nonempty: otherwise the closure starting from $S(P)$
would remain empty and could not reach $P$.

Suppose now that every $Mx$ is finite.  Each block $B_p^M$ ($p\in P$) has a maximal element rank; denote this maximum by
$\rho(p)$. Choose a total order $<$ on $P$ refining the strict
order of these maxima.  Fix an active $(p,q)$ and choose elements
$a,b$ in their respective blocks with ranks $\rho(p),\rho(q)$.
The pair $\langle a,b\rangle$ has a profile $r\in T(p,q)$ and
\[
 \rho(r)\ge\operatorname{rk}\langle a,b\rangle
 >\max\{\rho(p),\rho(q)\}.
\]
Thus $r$ lies after both $p$ and $q$ in the chosen order, and
\eqref{eq:weakorder} follows.
\end{proof}

\begin{lemma}[Initialization of profile blocks]
\label{lem:certinitial}
Given an ordered profile certificate $(P,T)$ with $P\ne\e$, one can
construct pairwise disjoint finite nonempty sets
$B_p\subseteq\HF$ ($p\in P$) such that
\begin{enumerate}[label=(\alph*)]
\item each $|B_p|\ge K:=|P|+1$;
\item every element of $B_p$ is either a \emph{source seed}
(a finite ordinal of cardinality at least three, occurring only
when $p\in S(P)$), or a pair $\langle a,b\rangle$ with
$a\in B_q$, $b\in B_r$ for an active $(q,r)$ and
$p\in T(q,r)$;
\item every assigned pair has a unique input-profile pair and
lies in one of its allowed target blocks; moreover, for each active
node not used during initialization, no pair in its current product
has yet been assigned to any block.
\end{enumerate}
\end{lemma}
\begin{proof}
We shall first populate all the blocks while keeping enough elements in each nonempty block to supply the targets of an active node. The choice $K=|P|+1$ ensures that two such input blocks yield at least $K^2$ pairs, more than the $K|T(p,q)|$ pairs needed for all the targets of a node $(p,q)$. We must also ensure that these pairs have not already been assigned.

To this end, initialize each $B_p$, for $p\in S(P)$, with $K$ distinct finite ordinals greater than two, choosing disjoint sets of ordinals for different source profiles.  Call these the \emph{source seeds}.  An
ordinal with at least three members is not a Kuratowski ordered pair,
which has at most two members.  Initially take $B_p=\varnothing$
for all other profiles.  Accessibility implies $S(P)\ne\varnothing$
when $P\ne\varnothing$.  By \eqref{eq:targets-not-sources}, source blocks retain only their seeds.

Throughout initialization, we maintain the following properties. At each stage
the blocks are disjoint, the only assigned elements which are not
source seeds are pairs generated by active nodes already processed,
and each such pair lies in an allowed target block of its generating
node.  Further, if $(p,q)$ has not been processed, then
\begin{equation}\label{eq:unused-node-fresh}
 (B_p\times B_q)\cap\bigcup_{r\in P}B_r=\varnothing.
\end{equation}
Indeed, an element in this intersection cannot be a source seed. It must therefore have been generated by some previously processed node
$(p',q')$.  Injectivity of Kuratowski pairing and disjointness of the
blocks force $(p,q)=(p',q')$, contrary to the assumption that $(p,q)$
was not processed.  This also establishes uniqueness of the
input-profile pair for any assigned pair, including when $p=q$.

Suppose that some block is still empty. To find a node that can be processed, let $C_0=S(P)$ and define
\[
 C_{k+1}=C_k\cup
 \bigcup_{\substack{(p,q)\in N(P)\\p,q\in C_k}}T(p,q).
\]
Accessibility gives $C_k=P$ for some $k$.  Choose the least $k$
for which $C_k$ contains the profile of a currently empty block.
Then $k>0$, every block indexed by $C_{k-1}$ is nonempty, and
there is an active node $(p,q)$ with $p,q\in C_{k-1}$ whose target
set contains an empty block.  This node has not been processed:
processing a node adds elements to \emph{all} its target blocks,
and no block is subsequently emptied.  Each of its input blocks
has at least $K$ elements, either from its initial seeding or from
a previous initialization step.

By \eqref{eq:unused-node-fresh}, all $|B_p\times B_q|\ge K^2$
pairs are globally unassigned.  Since
\[
 K|T(p,q)|\le K|P|=K(K-1)<K^2,
\]
we can select $K|T(p,q)|$ distinct pairs and distribute them in
disjoint packets of $K$ pairs, one packet to each block $B_r$ with
$r\in T(p,q)$.  The remaining pairs are left unassigned.  Each
packet consists of HF elements because its two components are HF.
This step preserves disjointness, admissible provenance, and the
freshness property for every other unprocessed node, by the same
injectivity argument.  At least one formerly empty block becomes nonempty. Repeat the step while some block remains empty. Thus at most $|P|$ steps suffice, and all resulting blocks are finite subsets of $\HF$ of size at least $K$.

For use in the two realization lemmas below, we record the resulting \emph{provenance invariant}.  It holds after initialization and persists under
any extension that inserts only previously unassigned ordered pairs
$\langle a,b\rangle$ into a block $B_r$ with
$a\in B_p$, $b\in B_q$, $(p,q)\in N(P)$ and $r\in T(p,q)$,
while preserving block disjointness.  Namely, every element that is
not a source seed has this permitted origin, and its two components
belong to uniquely determined input-profile blocks.  Consequently,
for every active $(p,q)$,
\begin{equation}\label{eq:provenance}
 (B_p\times B_q)\cap\bigcup_{r\in P}B_r
 \ \subseteq\ \bigcup_{r\in T(p,q)}B_r;
\end{equation}
for inactive $(p,q)$ the intersection on the left is empty.
Thus a pair outside its allowed target blocks is globally unassigned, ensuring freshness at every subsequent completion step.
\end{proof}

Having populated every block, we can now complete the products required by the certificate, without introducing any further profiles.

\begin{lemma}[Ordinary realization]\label{lem:certordinary}
Every ordered profile certificate for $\Phi$ yields an ordinary
model of $\Phi$.
\end{lemma}
\begin{proof}
The empty certificate yields the all-empty assignment. If $P\ne\e$, start with the finite blocks of
Lemma~\ref{lem:certinitial}.  Fix, for each active node $(p,q)$,
a designated target $t(p,q)\in T(p,q)$.  If there are active nodes,
visit them in a periodic schedule, so that every node is visited
infinitely often; if there are none, no completion is needed.

At a visit to $(p,q)$, form
\[
 D_{p,q}=(B_p\times B_q)\setminus
       \bigcup_{r\in T(p,q)}B_r
\]
and add $D_{p,q}$ to $B_{t(p,q)}$.  By the provenance invariant
\eqref{eq:provenance}, no element of $D_{p,q}$ belongs to
\emph{any} currently assigned block.  Thus the extension preserves
disjointness, and each new pair is assigned to a permitted target.
Its components remain in their uniquely determined input blocks,
so provenance is preserved inductively at every finite stage.
Each stage adds only finitely many HF pairs.

Write $B_p^{\infty}$ for the increasing union of the finite-stage
blocks.  These unions are pairwise disjoint.  If
$a\in B_p^{\infty}$ and $b\in B_q^{\infty}$ and $(p,q)$ is active,
both elements occur at some finite stage.  At a later visit to
$(p,q)$ their pair is either already in a permitted target block
or is added to one; hence
\begin{equation}\label{eq:ordinary-closure}
 B_p^{\infty}\times B_q^{\infty}
 \ \subseteq\ \bigcup_{r\in T(p,q)}B_r^{\infty}
 \qquad((p,q)\in N(P)).
\end{equation}
Conversely, provenance at finite stages ensures that every
non-seed element in the final blocks is a pair produced by an
active node and assigned to one of its permitted targets.  No pair
from an inactive node belongs to any final block.

Define
\begin{equation}\label{eq:profilemodel}
 Mx=\bigcup_{\substack{p\in P\\p(x)=1}}B_p^{\infty}.
\end{equation}
Condition~\ref{cert:boolean} establishes all Boolean equations
and disequalities.  Let $x=y\times z$ be any product equation.
If $a\in My$ and $b\in Mz$, their unique profiles $p,q$
satisfy $p(y)q(z)=1$, so $(p,q)$ is active.  By
\eqref{eq:ordinary-closure}, $\langle a,b\rangle\in B_r^{\infty}$
for some $r\in T(p,q)$; then \eqref{eq:targetprofile} gives
$r(x)=1$, and $\langle a,b\rangle\in Mx$.
Conversely, if $c\in Mx$, its profile $r$ satisfies $r(x)=1$.
It cannot be a source seed by \eqref{eq:sourceprofiles}.
Thus $c=\langle a,b\rangle$ arose from an active $(p,q)$ with
$r\in T(p,q)$, and \eqref{eq:targetprofile} gives
$p(y)=q(z)=1$.  Therefore $a\in My$ and $b\in Mz$.
This proves $Mx=My\times Mz$ for every product literal,
simultaneously.
\end{proof}

\begin{corollary}[Optimal rank bound for ordinary models]\label{cor:orderedrank}
Every ordinarily satisfiable $\BSTC$ formula admits a model $M$ with $Mx\subseteq\HF$ for each variable $x$, or equivalently $\operatorname{rk}(Mx)\le\omega$. In particular, every variable value can be taken to be at most countable. The uniform rank bound $\omega$ is optimal.
\end{corollary}
\begin{proof}
For a normalized conjunction, Lemma~\ref{lem:certnecessary} supplies an ordered profile certificate. If the certificate is empty, the all-empty assignment has the required rank bound. Otherwise, the initialization of Lemma~\ref{lem:certinitial} uses only HF elements, and the construction of Lemma~\ref{lem:certordinary} adds only Kuratowski pairs of existing HF elements. Thus every block at every finite stage is a subset of $\HF$, as is its union at the limit. Equation~\eqref{eq:profilemodel} then gives $Mx\subseteq\HF$ for each variable $x$. In either case, $Mx\subseteq\HF=\mathcal V_\omega$, equivalently $Mx\in\mathcal V_{\omega+1}=\Pow(\mathcal V_\omega)$. This establishes the rank bound and countability, since $\HF$ is countable. For an arbitrary satisfiable formula, choose the truth values of its atoms in a model, flatten the resulting conjunction of literals by Lemma~\ref{lem:flatten} (replacing negative inclusions by nonemptiness of differences), and restrict the constructed model to the original variables.

For optimality, recall that the conjunction \eqref{eq:ordinaryonly} forces $My$ to be infinite in every model, as shown in Example~\ref{ex:ordinaryonly}. Since each finite level of the cumulative hierarchy is finite, $\operatorname{rk}(My)\ge\omega$. The bound just proved is therefore optimal.
\end{proof}

To obtain a hereditarily finite model, we need to arrange the completion steps so that the input blocks of a completed node are not enlarged by subsequent steps. As the next lemma shows, the weak-order condition provides just such an arrangement.

\begin{lemma}[HF realization under a weak order]
\label{lem:certhf}
Every weakly ordered profile certificate for $\Phi$ yields a
hereditarily finite model of $\Phi$.
\end{lemma}
\begin{proof}
The empty-certificate case is immediate.  For nonempty $P$,
fix a total order $<$ witnessing weak ordering and initialize the finite HF blocks $B_p$ by
Lemma~\ref{lem:certinitial}.  Enumerate all active nodes as
$e_i=(p_i,q_i)$, $1\le i\le k$, in nondecreasing order of $h_i=\max_<\{p_i,q_i\}$, breaking ties arbitrarily.  Put $t_i=\max_<T(p_i,q_i)$.
By the weak-order requirement, $h_i<t_i$.

At step $i$ add every then missing pair from
\[
 (B_{p_i}\times B_{q_i})\setminus
     \bigcup_{r\in T(p_i,q_i)}B_r
\]
to $B_{t_i}$.  The provenance invariant
\eqref{eq:provenance} implies that each missing pair is unassigned
globally, so this operation preserves disjointness and provenance.
All blocks remain finite subsets of $\HF$ because only finitely
many pairs of existing HF elements are added at each step.

It remains to establish the \emph{stability} of each completed node.  Once $e_i$ has
been completed, any later step $j>i$ satisfies
\begin{equation}\label{eq:hf-stability}
 t_j>h_j\ge h_i=\max_<\{p_i,q_i\}.
\end{equation}
Consequently, step $j$ enlarges only the block $B_{t_j}$,
which is neither $B_{p_i}$ nor $B_{q_i}$.  Since initialization
precedes all completion steps, the two input blocks of $e_i$ never
change after step $i$.  At step $i$ every pair in their then-current
product is either already in an allowed target block by
\eqref{eq:provenance}, or is placed in $B_{t_i}$; this inclusion
therefore remains true in the final blocks:
\begin{equation}\label{eq:hf-closure}
 B_{p_i}\times B_{q_i}
 \ \subseteq\ \bigcup_{r\in T(p_i,q_i)}B_r
 \qquad(1\le i\le k).
\end{equation}
Only the designated target needs to follow both input profiles
(cf. Example~\ref{ex:multitarget}).

Define $Mx$ by \eqref{eq:profilemodel}, taking the final finite
blocks in place of $B_p^{\infty}$.  The Boolean and disequality
conditions hold by \ref{cert:boolean}.  For each product equation
$x=y\times z$, inclusion $My\times Mz\subseteq Mx$ follows
from \eqref{eq:hf-closure} and \eqref{eq:targetprofile}; the reverse
inclusion follows from the provenance invariant and the exclusion of
source seeds from product-variable values, exactly as in the proof
of Lemma~\ref{lem:certordinary}.  These arguments apply to all
product literals simultaneously.  Finally, each $Mx$ is a finite
subset of $\HF$ and therefore an element of $\HF$.
\end{proof}

\begin{theorem}[Finite-certificate characterization]
\label{thm:orderedcertificate}
A normalized $\BSTC$ conjunction is ordinarily
satisfiable if and only if it admits an ordered profile certificate.
It is finitely satisfiable if and only if it admits a weakly ordered
profile certificate, and the same characterization holds for
hereditarily finite satisfiability.  In particular, finite and
HF satisfiability coincide for this language.
\end{theorem}
\begin{proof}
Necessity is Lemma~\ref{lem:certnecessary}.  Sufficiency for
ordinary satisfiability is Lemma~\ref{lem:certordinary};
sufficiency for both finite and HF satisfiability is
Lemma~\ref{lem:certhf}.  The remaining implications are immediate.
The equivalence of finite and HF satisfiability extends to arbitrary
quantifier-free formulae: take the truth values of the atoms in a
finite model, replace negative inclusions, if used, by nonemptiness
constraints on differences, and apply Lemma~\ref{lem:flatten} to the
resulting conjunction of literals. An HF model of the normalized
conjunction restricts to an HF model of the original formula.
\end{proof}

For decision purposes, the characterization just obtained can be simplified. In fact, for a fixed set of profiles, we can take all the locally admissible targets without affecting the existence of a certificate or of a weak order.

\begin{proposition}[Canonical target relation]\label{prop:canonicaltargets}
For the fixed conjunction $\Phi$ and a set $P$ of nonzero profiles, define, for every $(p,q)\in N(P)$,
\begin{equation}\label{eq:canonicaltargets}
 T^*(p,q)=\{r\in P:r(x)=p(y)q(z)\text{ for every }(x=y\times z)\in\mathcal M\}.
\end{equation}
There is a target relation $T$ such that $(P,T)$ is an ordered profile certificate if and only if $(P,T^*)$ is one. The same holds for weakly ordered profile certificates.
\end{proposition}
\begin{proof}
If $(P,T)$ is a certificate, then $T(p,q)\subseteq T^*(p,q)$ by \eqref{eq:targetprofile}. Hence the canonical target sets are nonempty on active nodes and satisfy the product condition by definition. The Boolean and disequality conditions depend only on $P$. Accessibility is preserved under enlargement of target sets: starting from $S(P)$, every round of the closure under $T^*$ contains the corresponding round under $T$, by induction. Thus $(P,T^*)$ is accessible. If $<$ witnesses weak ordering for $(P,T)$, then for every active node
\[
 \max_<\{p,q\}<\max_<T^*(p,q),
\]
since the maximal canonical target is at least as large in $<$ as the maximal original target. The same order witnesses weak ordering for $(P,T^*)$. The converse in both cases follows by taking $T=T^*$.
\end{proof}

We show next that the weak-order condition can be checked directly on the incidence representation, without guessing a total order.

\begin{proposition}[Greedy recognition of weak orders]\label{prop:weakrecognition}
Given an ordered profile certificate $(P,T)$, the existence of a weak order can be decided deterministically, and a witnessing order constructed if one exists, in time $O(s+a+I)$ on an explicit incidence representation with indexed profiles and nodes, where $s=|P|$, $a=|N(P)|$, and $I=\sum_{(p,q)\in N(P)}|T(p,q)|$.
\end{proposition}
\begin{proof}
We construct the order from greatest to least. Let $A$ denote the set of profiles already placed (initially empty). A profile $p\in P\setminus A$ is eligible if every active node $(u,v)$ with $p\in\{u,v\}$ satisfies $T(u,v)\cap A\ne\e$. Repeatedly choose any eligible profile, place it below all profiles already placed, and add it to $A$. Accept if $A=P$, and reject if $A\ne P$ and no eligible profile exists.

To establish correctness, suppose first that the procedure accepts and let $<$ be its resulting order. For any active node $(p,q)$, when $h=\max_<\{p,q\}$ was placed, at least one of its targets was already in $A$, hence above $h$. This establishes \eqref{eq:weakorder}.

Conversely, suppose a witnessing order $\prec$ exists. At any stage with $A\ne P$, choose $p=\max_\prec(P\setminus A)$. For every active node with $p$ as an input, some target is strictly above both inputs in $\prec$, and in particular above $p$. Such a target belongs to $A$ by the choice of $p$. Thus $p$ is eligible, regardless of earlier greedy choices, and the procedure cannot reject.

Concerning the implementation and its complexity, call a node satisfied when at least one of its targets belongs to $A$. Construct the inverse target-incidence lists by scanning the explicit representation. Maintain for each profile a count of the unsatisfied nodes in which it occurs as an input, counting a diagonal node only once. Initially all nodes are unsatisfied; the profiles with count zero are placed in a queue. When a profile is selected, scan the list of nodes having it as a target. Whenever a node becomes satisfied for the first time, decrement the counts of its distinct inputs, enqueuing any unplaced profile whose count has just become zero. Counts are nonincreasing, so each profile enters the queue at most once. Each node changes status at most once and each target incidence is scanned at most once. Initialization and the subsequent operations take $O(s+a+I)$ time and space in the stated representation, with the empty case handled immediately. The bound counts unit-cost operations on indices; a standard bit encoding incurs at most polynomial overhead.
\end{proof}

\begin{corollary}[Characterization by profile sets]\label{cor:profilesonly}
A normalized $\BSTC$ conjunction is ordinarily satisfiable if and only if there is a set $P$ of nonzero profiles satisfying condition~\ref{cert:boolean} of Definition~\ref{def:cert}, with nonempty canonical target sets on all active nodes and accessibility under $T^*$. For finite and HF satisfiability, the same conditions apply with the additional requirement that the greedy procedure of Proposition~\ref{prop:weakrecognition} accepts $(P,T^*)$.
\end{corollary}
\begin{proof}
Combine Theorem~\ref{thm:orderedcertificate} with Propositions~\ref{prop:canonicaltargets} and~\ref{prop:weakrecognition}.
\end{proof}

\begin{proposition}[\NEXPTIME upper bound]\label{prop:orderedupper}
Ordinary, finite, and hereditarily finite satisfiability for
$\BSTC$ belong to \NEXPTIME, for normalized
conjunctions and for arbitrary quantifier-free formulae.
\end{proposition}
\begin{proof}
Let $\Phi$ be a normalized conjunction with $m$ variables. By Corollary~\ref{cor:profilesonly}, it is enough to guess a subset $P$ of the at most $2^m-1$ nonzero profiles. Use a characteristic vector of $2^m-1$ bits indexed by these profiles, and scan it to obtain an explicit list of the selected profiles, in time $O(m2^m)$ for $m\ge 1$ (constant time for $m=0$).

Verify condition~\ref{cert:boolean} and compute $N(P)$ and the canonical target relation \eqref{eq:canonicaltargets} deterministically. Reject if any active node has an empty target set or if the accessibility test by iterative closure fails. For finite and HF satisfiability, apply in addition the greedy test of Proposition~\ref{prop:weakrecognition}. The same corollary establishes the correctness of these tests: neither the target relation nor a witnessing order needs to be guessed, and we need not construct the HF model itself, whose elements may be far more numerous than the certificate.

Concerning complexity, there are at most $|P|^2$ nodes and $|P|^3$ target incidences. The verification of the local conditions and the construction of the canonical target relation thus take time polynomial in $|\Phi|$ and $|P|+1$. The accessibility test requires at most $|P|$ strict enlargement rounds, and the greedy test is polynomial in the explicit incidence representation. Taking into account the cost of representing the profiles, the entire test takes at most $2^{O(m)}\operatorname{poly}(|\Phi|)$ time.

For an arbitrary quantifier-free formula, guess a truth value
for each atomic comparison, check the resulting propositional
skeleton, and flatten the corresponding positive and negative
set-term equalities using fresh variables.  Negated equality
becomes a disequality, and negative inclusions become
nonemptiness constraints on set differences.  This produces,
in polynomial time, a normalized conjunction preserving each
of the three semantics.  The stated upper bound follows.
\end{proof}

\begin{corollary}[Decidability]\label{cor:decidability}
Ordinary, finite, and hereditarily finite satisfiability are decidable
for both $\BSTuC$ and $\BSTC$,
considered as separate languages.
\end{corollary}
\begin{proof}
The result follows from Proposition~\ref{prop:orderedupper} and the \NEXPTIME upper bounds established in Cantone and Ursino~\cite{CantoneUrsino2024}.
\end{proof}

\begin{theorem}[Ordered Cartesian product]\label{thm:orderedmain}
Ordinary, finite, and hereditarily finite satisfiability for
$\BSTC$ are \NEXPTIME-complete, already for
normalized conjunctions of atomic set constraints and
set disequalities.
\end{theorem}
\begin{proof}
Combine Proposition~\ref{prop:orderedhard} with
Proposition~\ref{prop:orderedupper}.
\end{proof}

\begin{corollary}\label{cor:notnp}
None of the six satisfiability problems for $\BSTuC$ and $\BSTC$, considered as separate languages, belongs to
$\NP$.
\end{corollary}
\begin{proof}
The nondeterministic time hierarchy theorem implies
$\NP\subsetneq\NEXPTIME$. If any of these \NEXPTIME-hard
problems were in \NP, closure of \NP under polynomial-time many-one
reductions would imply $\NEXPTIME\subseteq\NP$, a contradiction.
\end{proof}

\section{Discussion and conclusions}
We have established that ordinary, finite, and hereditarily finite satisfiability are \NEXPTIME-complete for both alternative Cartesian-product extensions of Boolean set theory. For the unordered language, our reduction matches the upper bounds of Cantone and Ursino~\cite{CantoneUrsino2024}; for the ordered language, we establish both lower and upper bounds.

In the lower-bound construction, address generation and adjacency testing are treated separately, with the source and destination sets determining the orientation of unordered successor pairs. For the ordered product, our characterization reduces the nondeterministic choice for a normalized conjunction to the set of profiles (Corollary~\ref{cor:profilesonly}); the weak-order condition of Definition~\ref{def:cert} distinguishes ordinary from finite and HF satisfiability. Multivalued targets of arbitrarily large cardinality can be necessary (Proposition~\ref{prop:targetmultiplicity}). Ordinary satisfiability admits the optimal uniform rank bound $\omega$ on variable values (Corollary~\ref{cor:orderedrank}). Although the certificates are exponentially bounded, the realization may produce much larger finite models: the certificate bound is not an explicit small-model bound. Nor does it establish a deterministic single-exponential algorithm. The characterization and its realization proofs do not rely on a reduction to unordered-product satisfiability.

The results concern the two languages separately, not the combined language with both products. Extensions with membership predicates or cardinality comparisons likewise require separate analysis.

\appendix
\section{A worked two-bit grid and its tiling constraints}
\label{app:twobit}
We illustrate the reduction on the smallest nontrivial grid, giving the HF branching witnesses, the successor relations, and one choice of tile constraints. In particular, we show how unordered pairs can represent directed adjacencies.

\paragraph{The two branching stages.}
Take $n=1$, so the square has side $2$, and let $k=2$.
For a particular HF model of $\Gamma_2$, put
\[
 L_0=\{\e\},\qquad
 T_0^0=\{3\},\quad T_0^1=\{4\},\qquad
 T_1^0=\{5\},\quad T_1^1=\{6\},
\]
where the numbers denote finite von Neumann ordinals. The first
stage produces the two distinct elements
\[
 v_0=\{\e,3\},\qquad v_1=\{\e,4\},\qquad
 L_1=\{v_0,v_1\}.
\]
The second stage produces four final elements:
\[
\begin{array}{ll}
 g_{00}=\{v_0,5\},&g_{10}=\{v_1,5\},\\[2pt]
 g_{01}=\{v_0,6\},&g_{11}=\{v_1,6\}.
\end{array}
\]
The last-level set is $G=L_2=\{g_{00},g_{10},g_{01},g_{11}\}$.
All four elements are distinct and hereditarily finite. The bit
sets are
\[
 B_1=\{g_{10},g_{11}\},\qquad
 B_2=\{g_{01},g_{11}\}.
\]
The first bit is the horizontal coordinate $x$, and the second
is the vertical coordinate $y$. Thus the four singleton address
cells have the following interpretation:
\[
\begin{array}{c|c|c|c}
 \text{element}&(\beta_1,\beta_2)&(x,y)&\text{cell}\\\hline
 g_{00}&(0,0)&(0,0)&\mathcal{R}_{00}\\
 g_{10}&(1,0)&(1,0)&\mathcal{R}_{10}\\
 g_{01}&(0,1)&(0,1)&\mathcal{R}_{01}\\
 g_{11}&(1,1)&(1,1)&\mathcal{R}_{11}
\end{array}
\]

\paragraph{Successor relations.}
In the notation of Section~\ref{sec:grid}, the horizontal and
vertical pivots are $i=1$ and $i=2$, with
\[
 S_1=\bar B_1,\quad D_1=B_1,\quad J_1=\{2\},\qquad
 S_2=\bar B_2,\quad D_2=B_2,\quad J_2=\{1\}.
\]
Consequently, the unordered successor-pair sets are
\begin{align*}
 R_1 &=(\bar B_1\otimes B_1)\cap Q_{\{2\}}
      =\bigl\{\{g_{00},g_{10}\},\{g_{01},g_{11}\}\bigr\},\\
 R_2 &=(\bar B_2\otimes B_2)\cap Q_{\{1\}}
      =\bigl\{\{g_{00},g_{01}\},\{g_{10},g_{11}\}\bigr\}.
\end{align*}
No diagonal pair, such as $\{g_{00},g_{11}\}$, belongs to either
relation: both coordinates would change. The sets $S_i$ and $D_i$
fix the unique orientation of each unordered pair, from bit $0$
to bit $1$ in the pivot position.

\paragraph{An explicit tiling.}
Choose $D=\{A,B\}$,
$H=V=\{(A,B),(B,A)\}$, and $d_\star=A$.
The checkerboard assignment
\[
\begin{array}{c|cc}
 &x=0&x=1\\\hline
 y=1&B&A\\
 y=0&A&B
\end{array}
\qquad\text{corresponds to}\qquad
 C_A=\{g_{00},g_{11}\},\quad
 C_B=\{g_{10},g_{01}\}.
\]
It meets the origin condition, and each directed successor pair
has one of the permitted color combinations $(A,B)$ or $(B,A)$.
For instance, the horizontal successor from $(0,0)$ to $(1,0)$
is represented by $\{g_{00},g_{10}\}$, where the source is in
$S_1\cap C_A$ and the destination is in $D_1\cap C_B$.
The forbidden combination $(A,A)$ therefore yields the clause
\[
 \bigl((S_1\cap C_A)\otimes(D_1\cap C_A)\bigr)
 \cap Q_{\{2\}}=\e,
\]
and the checkerboard interpretation satisfies it. By contrast,
coloring all four cells $A$ would violate this clause, since
$\{g_{00},g_{10}\}$ would lie in the forbidden intersection.
Thus the example illustrates both directions of the reduction.

\section{A carry pivot beyond the smallest grid}
\label{app:carry}
For $n=2$ the addresses have four bits $(b_1,b_2,b_3,b_4)$,
with the first two bits encoding $x$ in least-significant-bit-first
order, and the last two encoding $y$. Consider the horizontal
successor $x=1\longrightarrow x=2$ at fixed $y$:
\[
 (1,0,b_3,b_4)\longrightarrow(0,1,b_3,b_4).
\]
Here the pivot is $i=2$ and the lower bit changes from 1 to 0 and generates a carry into the pivot position.
The source and destination conditions become
\[
 S_2=\bar B_2\cap B_1,\qquad
 D_2=B_2\cap\bar B_1,\qquad
 J_2=\{3,4\}.
\]
Thus $Q_{J_2}$ checks exactly the two vertical bits, leaving
the lower horizontal bit to be handled by the carry conditions.
For $i=1$, by contrast, the pivot simply flips the lowest
horizontal bit and preserves the other three bits. Hence the reduction distinguishes the relations $R_i$ according to the carry pivot, though their union consists precisely of the unordered pairs of grid elements whose coordinates are horizontally or vertically adjacent.

\section*{Acknowledgments}

\noindent\textit{Declaration of generative AI use.}
OpenAI's ChatGPT was used in the preparation of this
manuscript to assist with checking and refining
mathematical arguments, bibliographic searches,
and linguistic and editorial revision.
The authors take full responsibility for the accuracy
of the mathematical results and references, and for
the final version of the manuscript.

\end{document}